\documentclass[11pt]{article}
\usepackage{amssymb}
\usepackage{amsmath}
\usepackage[all]{xy}

\title{\bf On $U$-Dominant Dimension
}
\author{Zhaoyong Huang\thanks{\small \it E-mail address: huangzy@nju.edu.cn}\\
{\small \it Department of Mathematics, Nanjing University,}\\
{\small \it Nanjing 210093, The People's Republic of China}\\}
\usepackage{amssymb}
\date{}
\begin{document}
\baselineskip=18pt \maketitle

\begin{abstract}
Let $\Lambda$ and $\Gamma$ be artin algebras and $_{\Lambda}U
_{\Gamma}$ a faithfully balanced selforthogonal bimodule. We show
that the $U$-dominant dimensions of $_{\Lambda}U$ and $U
_{\Gamma}$ are identical. As applications to the results obtained,
we give some characterizations of double dual functors (with
respect to $_{\Lambda}U _{\Gamma}$) preserving monomorphisms and
being left exact respectively.
\end{abstract}

\vspace{1cm}

\centerline{\bf 1. Introduction}

\vspace{0.2cm}

For a ring $\Lambda$, we use mod $\Lambda$ (resp. mod $\Lambda
^{op}$) to denote the category of finitely generated left
$\Lambda$-modules (resp. right $\Lambda$-modules).

\vspace{0.2cm}

{\bf Definition 1.1} Let $\Lambda$ and $\Gamma$ be rings. A
bimodule $_{\Lambda}T _{\Gamma}$ is called a faithfully balanced
selforthogonal bimodule if it satisfies the following conditions:

(1) $_{\Lambda}T \in$mod $\Lambda$ and $T_{\Gamma} \in$mod $\Gamma
^{op}$.

(2) The natural maps $\Lambda \rightarrow {\rm End}(T_{\Gamma})$
and $\Gamma \rightarrow {\rm End}(_{\Lambda}T)^{op}$ are
isomorphisms.

(3) ${\rm Ext}_{\Lambda}^{i}(_{\Lambda}T , {_{\Lambda}T})=0$ and
${\rm Ext}_{\Gamma}^{i}(T_{\Gamma}, T_{\Gamma})=0$ for any $i \geq
1$.

\vspace{0.2cm}

{\bf Definition 1.2} Let $U$ be in mod $\Lambda$ (resp. mod
$\Gamma ^{op}$) and $n$ a non-negative integer. For a module $M$
in mod $\Lambda$ (resp. mod $\Gamma ^{op}$),

(1)$^{[8]}$ $M$ is said to have $U$-dominant dimension greater
than or equal to $n$, written $U$-dom.dim$(_{\Lambda}M)$ (resp.
$U$-dom.dim$(M_{\Gamma}))\geq n$, if each of the first $n$ terms
in a minimal injective resolution of $M$ is cogenerated by
$_{\Lambda}U$ (resp. $U _{\Gamma}$), that is, each of these terms
can be embedded into a direct product of copies of $_{\Lambda}U$
(resp. $U _{\Gamma}$).

(2)$^{[10]}$ $M$ is said to have dominant dimension greater than
or equal to $n$, written \linebreak dom.dim$(_{\Lambda}M)$ (resp.
dom.dim$(M_{\Gamma}))\geq n$, if each of the first $n$ terms in a
minimal injective resolution of $M$ is $\Lambda$-projective (resp.
$\Gamma ^{op}$-projective).

\vspace{0.2cm}

Assume that $\Lambda$ is an artin algebra. By [4] Theorem 3.3,
$\Lambda ^I$ and each of its direct summands are projective for
any index set $I$. So, when $_{\Lambda}U ={_{\Lambda}\Lambda}$
(resp. $U _{\Gamma}= \Lambda _{\Lambda}$), the notion of
$U$-dominant dimension coincides with that of (classical) dominant
dimension. Tachikawa in [10] showed that if $\Lambda$ is a left
and right artinian ring then the dominant dimensions of
$_{\Lambda}\Lambda$ and $\Lambda _{\Lambda}$ are identical. Kato
in [8] characterized the modules with $U$-dominant dimension
greater than or equal to one. Colby and Fuller in [5] gave some
equivalent conditions of dom.dim$(_{\Lambda}\Lambda)\geq$1 (or 2)
in terms of the properties of double dual functors (with respect
to ${_{\Lambda}\Lambda} _{\Lambda}$).

The results mentioned above motivate our interests in establishing
the identity of $U$-dominant dimensions of $_{\Lambda}U$ and $U
_{\Gamma}$ and characterizing the properties of modules with a
given $U$-dominant dimension. Our characterizations will lead a
better comprehension about $U$-dominant dimension and the theory
of selforthogonal bimodules.

Throughout this paper, $\Lambda$ and $\Gamma$ are artin algebras
and $_{\Lambda}U _{\Gamma}$ is a faithfully balanced
selforthogonal bimodule. The main result in this paper is the
following

\vspace{0.2cm}

{\bf Theorem 1.3} $U$-dom.dim$(_{\Lambda}U)=U$-dom.dim$(U
_{\Gamma})$.

\vspace{0.2cm}

Put $_{\Lambda}U_{\Gamma}={_{\Lambda}\Lambda _{\Lambda}}$, we
immediately get the following result, which is due to Tachikawa
(see [10]).

\vspace{0.2cm}

{\bf Corollary 1.4} dom.dim$(_{\Lambda}\Lambda)=$dom.dim$(\Lambda
_{\Lambda})$.

\vspace{0.2cm}

Let $M$ be in mod $\Lambda$ (resp. mod $\Gamma ^{op}$) and $G(M)$
the subcategory of mod $\Lambda$ (resp. mod $\Gamma ^{op}$)
consisting of all submodules of the modules generated by $M$. $M$
is called a QF-3 module if $G(M)$ has a cogenerator which is a
direct summand of every other cogenerator$^{[11]}$. By [11]
Proposition 2.2 we have that a finitely cogenerated
$\Lambda$-module (resp. $\Gamma ^{op}$-module) $M$ is a QF-3
module if and only if $M$ cogenerates its injective envelope. So
by Theorem 1.3 we have

\vspace{0.2cm}

{\bf Corollary 1.5} $_{\Lambda}U$ {\it is} QF-3 {\it if and only
if} $U _{\Gamma}$ {\it is} QF-3.

\vspace{0.2cm}

We shall prove our main result in Section 2. As applications to
the results obtained in Section 2, we give in Section 3 some
characterizations of double dual functors (with respect to
$_{\Lambda}U _{\Gamma}$) preserving monomorphisms and being left
exact respectively.

\vspace{0.5cm}

\centerline{\bf 2. The proof of main result}

\vspace{0.2cm}

Let $E_{0}$ be the injective envelope of $_{\Lambda} U$. Then
$E_{0}$ defines a torsion theory in mod $\Lambda$. The torsion
class $\cal T$ is the subcategory of mod $\Lambda$ consisting of
the modules $X$ satisfying Hom$_{\Lambda}(X,E_{0})=0$, and the
torsionfree class $\cal F$ is the subcategory of mod $\Lambda$
consisting of the modules $Y$ cogenerated by $E_{0}$
(equivalently, $Y$ can be embedded in $E_{0}^{I}$ for some index
set $I$). A module in mod $\Lambda$ is called torsion (resp.
torsionfree) if it is in $\cal T$ (resp. $\cal F$). The injective
envelope $E_0^{'}$ of $U _{\Gamma}$ also defines a torsion theory
in mod $\Gamma ^{op}$ and we may give in mod $\Gamma ^{op}$ the
corresponding notions as above. Let $X$ be in mod $\Lambda$ (resp.
mod $\Gamma ^{op}$) and $t(X)$ the torsion submodule, that is,
$t(X)$ is the submodule $X$ such that
Hom$_{\Lambda}(t(X),E_{0})=0$ (resp.
Hom$_{\Gamma}(t(X),E_{0}^{'})=0$) and $E_0$ (resp. $E_0^{'}$)
cogenerates $X/t(X)$ (c.f. [7]).

Let $A$ be in mod $\Lambda$ (resp. mod $\Gamma ^{op}$). We call
Hom$_{\Lambda}(_{\Lambda}A, {_{\Lambda}U _{\Gamma}})$ (resp.
Hom$_{\Gamma}(A_{\Gamma}, {_{\Lambda}U_{\Gamma}})$) the dual
module of $A$ with respect to $_{\Lambda}U_{\Gamma}$, and denote
either of these modules by $A^*$. For a homomorphism $f$ between
$\Lambda$-modules (resp. $\Gamma ^{op}$-modules), we put $f^*={\rm
Hom}(f, {_{\Lambda}U} _{\Gamma})$. Let $\sigma _{A}: A \rightarrow
A^{**}$ via $\sigma _{A}(x)(f)=f(x)$ for any $x \in A$ and $f \in
A^*$ be the canonical evaluation homomorphism. $A$ is called
$U$-torsionless (resp. $U$-reflexive) if $\sigma _{A}$ is a
monomorphism (resp. an isomorphism).

\vspace{0.2cm}

{\bf Lemma 2.1} {\it For a module} $X$ {\it in} mod $\Lambda$
({\it resp.} mod $\Gamma ^{op}$), $t(X)=$Ker$\sigma _{X}$ {\it if
and only if} Hom$_{\Lambda}($Ker$\sigma _{X}, E_{0})=0$ ({\it
resp.} Hom$_{\Gamma}($Ker$\sigma _{X}, E_{0}^{'})=0$).

\vspace{0.2cm}

{\it Proof.} The necessity is trivial. Now we prove the
sufficiency.

We have the following commutative diagram with the upper row
exact:

$$\xymatrix{0 \ar[r] & t(X) \ar[r] & X \ar[r]^{\pi} \ar[d]^{\sigma _X} &
X/t(X) \ar[r] \ar[d]^{\sigma _{X/t(X)}} & 0 \\
& & X^{**} \ar[r]^{\pi ^{**}} & [X/t(X)]^{**} &}
$$
Since Hom$_{\Lambda}(t(X), E_{0})=0$, $[t(X)]^*=0$ and $\pi ^*$ is
an isomorphism. So $\pi ^{**}$ is also an isomorphism and hence
$t(X)\subset$Ker$\sigma _{X}$. On the other hand,
Hom$_{\Lambda}($Ker$\sigma _{X}, E_{0})=0$ by assumption, which
implies that Ker$\sigma _{X}$ is a torsion module and contained in
$X$. So we conclude that Ker$\sigma _{X} \subset t(X)$ and
Ker$\sigma _{X}=t(X)$.  $\blacksquare$

\vspace{0.2cm}

{\it Remark.} From the above proof we always have
$t(X)\subset$Ker$\sigma _{X}$.

\vspace{0.2cm}

Suppose that $A\in$mod $\Lambda$ (resp. mod $\Gamma ^{op}$) and
$P_{1} \buildrel {f} \over \longrightarrow P_{0} \to A \to 0$ is a
(minimal) projective resolution of $A$. Then we have an exact
sequence $0 \to A^* \to P_{0}^* \buildrel {f^*} \over
\longrightarrow P_{1}^* \to {\rm Coker}f^* \to 0$. We call
Coker$f^*$ the transpose (with respect to $_{\Lambda}U _{\Gamma}$)
of $A$, and denote it by Tr$_{U}A$.

\vspace{0.2cm}

{\bf Proposition 2.2} {\it The following statements are
equivalent.}

(1) $t(X)=$Ker$\sigma _{X}$ {\it for every} $X \in$mod $\Lambda$.

(2) $f^{**}$ {\it is monic for every monomorphism} $f:A \to B$
{\it in} mod $\Lambda$.

$(1)^{op}$ $t(Y)=$Ker$\sigma _{Y}$ {\it for every} $Y \in$mod
$\Gamma ^{op}$.

$(2)^{op}$ $g^{**}$ {\it is monic for every monomorphism} $g:C \to
D$ {\it in} mod $\Gamma ^{op}$.

\vspace{0.2cm}

{\it Proof.} By symmetry, it suffices to prove the implications of
$(1)\Rightarrow (2)^{op} \Rightarrow (1)^{op}$.

$(1) \Rightarrow (2)^{op}$ Let $g:C\to D$ be monic in mod $\Gamma
^{op}$. Set $X=$Coker$g$. We have that Ker$\sigma _{{\rm
Tr}_{U}X}\cong$Ext$_{\Lambda}^{1}(X, U)$ and Tr$_{U}X \in$mod
$\Lambda$ by [6] Lemma 2.1. By (1) and Lemma 2.1,
Hom$_{\Gamma}($Ext$_{\Lambda}^{1}(X, U), E_{0})=0$. Since
Coker$g^*$ can be imbedded in Ext$_{\Lambda}^{1}(X, U)$,
Hom$_{\Gamma}$(Coker$g^*,$\linebreak $E_{0})=0$. But
(Coker$g^*)^*\subset$ Hom$_{\Gamma}$(Coker$g^*, E_{0})$, so
(Coker$g^*)^*=0$ and hence Ker$g^{**}\cong$\linebreak
(Coker$g^*)^*=0$, which implies that $g^{**}$ is monic.

$(2)^{op}\Rightarrow (1)^{op}$ Let $Y$ be in mod $\Gamma ^{op}$
and X any submodule of Ker$\sigma _{Y}$ and $f_{1}:
X\to$Ker$\sigma _{Y}$ the inclusion. Assume that $f$ is the
composition: $X \buildrel {f_{1}} \over \longrightarrow$Ker$\sigma
_{Y} \to Y$. Then $\sigma _{Y}f=0$ and $f^*\sigma _{Y}^*=(\sigma
_{Y}f)^*=0$. But $\sigma _{Y}^*$ is epic by [1] Proposition 20.14,
so $f^*=0$ and $f^{**}=0$. By $(2)^{op}$, $f^{**}$ is monic, so
$X^{**}=0$ and $X^{***}=0$. Since $X^*$ is isomorphic to a
submodule of $X^{***}$ by [1] Proposition 20.14, $X^*=0$.

We claim that Hom$_{\Gamma}($Ker$\sigma _{Y}, E_{0}^{'})=0$.
Otherwise, there exists $0\neq \alpha \in$
Hom$_{\Gamma}($Ker$\sigma _{Y}, E_{0}^{'})$. Then Im$\alpha
\bigcap U _{\Gamma} \neq 0$ since $ U _{\Gamma}$ is an essential
submodule of $E_{0}^{'}$. So $\alpha ^{-1}$(Im$\alpha \bigcap U
_{\Gamma})$ is a non-zero submodule of Ker$\sigma _{Y}$ and there
exists a non-zero map $\alpha ^{-1}$(Im$\alpha \bigcap U
_{\Gamma})\to U _{\Gamma}$, which implies that $(\alpha
^{-1}$(Im$\alpha \bigcap U _{\Gamma}))^*\neq 0$, a contradiction
with the former argument. Hence we conclude that $t(Y)=$Ker$\sigma
_{Y}$ by Lemma 2.1. $\blacksquare$

\vspace{0.2cm}

Let $A$ be a $\Lambda$-module (resp. a $\Gamma ^{op}$-module). We
denote either of Hom$_{\Lambda}(_{\Lambda}U _{\Gamma},
{_{\Lambda}A})$ and Hom$_{\Gamma}(_{\Lambda}U _{\Gamma},
A_{\Gamma})$ by $^*A$, and the left (resp. right) flat dimension
of $A$ by {\it l.}fd$_{\Lambda}(A)$ (resp. {\it
r.}fd$_{\Gamma}(A)$). We give a remark as follows. For an artin
algebra $R$ and a left (resp. right) $R$-module $A$, we have that
the left (resp. right) flat dimension of $A$ and its left (resp.
right) projective dimension are identical; especially, $A$ is left
(resp. right) flat if and only if it is left (resp. right)
projective.

\vspace{0.2cm}

{\bf Lemma 2.3} {\it Let} $_{\Lambda}E$ ({\it resp.} $E_{\Gamma
}$) {\it be injective and} $n$ {\it a non-negative integer. Then}
{\it l.}fd$_{\Gamma}({^*E})$ ({\it resp.} {\it
r.}fd$_{\Lambda}({^*E}))\leq n$ {\it if and only if}
Hom$_{\Lambda}($Ext$_{\Gamma}^{n+1}(A, U), E)$ ({\it resp.}
Hom$_{\Gamma}($Ext$_{\Lambda}^{n+1}(A, U), E))$ $=0$ {\it for any}
$A \in$mod $\Gamma ^{op}$ ({\it resp.} mod $\Lambda$).

\vspace{0.2cm}

{\it Proof.} It is trivial by [3] Chapter VI, Proposition 5.3.
$\blacksquare$

\vspace{0.2cm}

{\bf Proposition 2.4} {\it The following statements are
equivalent.}

(1) ${^*E}_{0}$ {\it is flat}.

(2) {\it There is an injective} $\Lambda$-{\it module} $E$ {\it
such that} ${^*E}$ {\it is flat and} $E$ {\it cogenerates}
$E_{0}$.

(3) $t(X)=$Ker$\sigma _{X}$ {\it for any} $X \in$mod $\Lambda$.

\vspace{0.2cm}

{\it Proof.} $(1) \Rightarrow (2)$ It is trivial.

$(2) \Rightarrow (3)$ Let $X \in$mod $\Lambda$. Since Ker$\sigma
_{X}\cong$Ext$^{1}_{\Gamma}($Tr$_{U}X, U)$ with Tr$_{U}X\in$mod
$\Gamma ^{op}$ by [6] Lemma 2.1. By (2) and Lemma 2.3,
Hom$_{\Lambda}($Ext$^{1}_{\Gamma}($Tr$_{U}X, U), E)=0$.

Since $E$ cogenerates $E_{0}$, there is an exact sequence $0 \to
E_{0} \to E^{I}$ for some index set $I$. So
Hom$_{\Lambda}($Ext$^{1}_{\Gamma}$(Tr$_{U}X, U), E_{0}) \subset$
Hom$_{\Lambda}($Ext$^{1}_{\Gamma}($Tr$_{U}X, U), E^{I}) \cong$
[Hom$_{\Lambda}($Ext$^{1}_{\Gamma}($Tr$_{U}X, U), E)]^{I}=0$ and
Hom$_{\Lambda}($Ext$^{1}_{\Gamma}($Tr$_{U}X, U), E_{0})=0$. By
Lemma 2.1, $t(X)=$Ker$\sigma _{X}$.

$(3) \Rightarrow (1)$ Let $N \in$mod $\Gamma ^{op}$. Since
Ker$\sigma _{{\rm Tr}_{U}N}\cong$Ext$^{1}_{\Gamma}(N, U)$ with
Tr$_{U}N\in$mod $\Lambda$ by [6] Lemma 2.1, By (3) and Lemma 2.1
we have Hom$_{\Lambda}($Ext$^{1}_{\Gamma}(N, U), E_{0})\cong$
Hom$_{\Lambda}($Ker$\sigma _{{\rm Tr}_{U}N},
E_{0})$\linebreak
$=0$, and so ${^*E}_{0}$ is flat by Lemma 2.3.
$\blacksquare$

\vspace{0.2cm}

Dually, we have the following

\vspace{0.2cm}

{\bf Proposition 2.4}$'$ {\it The following statements are
equivalent.}

(1) ${^*E}_{0}'$ {\it is flat}.

(2) {\it There is an injective} $\Gamma ^{op}$-{\it module} $E'$
{\it such that} ${^*E}'$ {\it is flat and} $E'$ {\it cogenerates}
$E_{0}'$.

(3) $t(Y)=$Ker$\sigma _{Y}$ {\it for any} $Y \in$mod $\Gamma
^{op}$.

\vspace{0.2cm}

{\bf Corollary 2.5} ${^*E}_{0}$ {\it is flat if and only if}
${^*E}_{0}^{'}$ {\it is flat}.

\vspace{0.2cm}

{\it Proof.} By Propositions 2.2, 2.4 and 2.4$'$. $\blacksquare$

\vspace{0.2cm}

Let $A \in$ mod $\Lambda$ (resp. mod $\Gamma ^{op}$) and $i$ a
non-negative integer. We say that the grade of $A$ with respect to
$_{\Lambda}U_{\Gamma}$, written grade$_U A$, is greater than or
equal to $i$ if Ext$_{\Lambda}^{j}(A, U)=0$ (resp.
Ext$_{\Gamma}^{j}(A, U)=0$) for any $0 \leq j < i$.

\vspace{0.2cm}

{\bf Lemma 2.6} {\it Let} $X$ be in mod $\Gamma ^{op}$ {\it and}
$n$ {\it a non-negative integer. If} grade$_{U}X \geq n$ {\it and}
grade$_{U}$Ext$_{\Gamma}^{n}(X, U)\geq n+1$, {\it then}
Ext$_{\Gamma}^{n}(X, U)=0$.

\vspace{0.2cm}

{\it Proof.} Since $X^*$ is $U$-torsionless, $X^{**}=0$ if and
only if $X^*=0$. Then the case $n=0$ follows.

Now let $n \geq 1$ and
$$\cdots \to P_{n} \to \cdots \to P_{1} \to P_{0} \to X \to 0$$
be a projective resolution of $X$ in mod $\Gamma ^{op}$. Put
$X_{n}=$Coker$(P_{n+1}\to P_{n})$. Then we have an exact sequence
$$0 \to P_{0}^* \to \cdots \to P_{n-1}^*
\buildrel {f} \over \longrightarrow X_{n}^* \to {\rm
Ext}_{\Gamma}^{n}(X, U) \to 0$$ in mod $\Lambda$ with each
$P_{i}^* \in$add$_{\Lambda}U$. Since
grade$_{U}$Ext$_{\Gamma}^{n}(X, U)\geq n+1$,
Ext$_{\Lambda}^{i}$(Ext$_{\Gamma}^{n}(X, U), U)=0$ for any $0 \leq
i \leq n$. So Ext$_{\Lambda}^{i}$(Ext$_{\Gamma}^{n}(X, U),
P_{j}^*)=0$ for any $0 \leq i \leq n$ and $0 \leq j \leq n-1$, and
hence Ext$_{\Lambda}^{1}$(Ext$_{\Gamma}^{n}(X, U), {\rm
Im}f)\cong$ Ext$_{\Lambda}^{n}$(Ext$_{\Gamma}^{n}(X, U),
P_{0}^*)=0$, which implies that we have an exact sequence
Hom$_{\Lambda}$(Ext$_{\Gamma}^{n}(X, U), X_{n}^*) \to$
Hom$_{\Lambda}$(Ext$_{\Gamma}^{n}(X, U)$, Ext$_{\Gamma}^{n}(X, U))
\to 0$. Notice that $X_{n}^*$ is $U$-torsionless and
Hom$_{\Lambda}$(Ext$_{\Gamma}^{n}(X, U), U)=0$. So
Hom$_{\Lambda}$(Ext$_{\Gamma}^{n}(X, U), X_{n}^*)=0$ and
Hom$_{\Lambda}$(Ext$_{\Gamma}^{n}(X, U)$, Ext$_{\Gamma}^{n}(X, U))
=0$, which implies that Ext$_{\Gamma}^{n}(X, U)=0$. $\blacksquare$

\vspace{0.2cm}

{\it Remark.} We point out that all of the above results (from 2.1
to 2.6) in this section also hold in the case $\Lambda$ and
$\Gamma$ are left and right noetherian rings.

\vspace{0.2cm}

For a module $T$ in mod $\Lambda$ (resp. mod $\Gamma ^{op}$), we
use add$_{\Lambda}T$ (resp. add$T_{\Gamma}$) to denote the
subcategory of mod $\Lambda$ (resp. mod $\Gamma ^{op}$) consisting
of all modules isomorphic to direct summands of finite direct sums
of copies of $_{\Lambda}T$ (resp. $T_{\Gamma}$). Let $A$ be in mod
$\Lambda$. If there is an exact sequence $\cdots \to U _{n} \to
\cdots \to U _{1} \to U _{0} \to A \to 0$ in mod $\Lambda$ with
each $U _{i}\in$add$_{\Lambda}U$ for any $i \geq 0$, then we
define $U$-resol.dim$_{\Lambda}(A)=$inf$\{ n|$ there is an exact
sequence $0 \to U _{n} \to \cdots \to U _{1} \to U _{0} \to A \to
0$ in mod $\Lambda$ with each $U _{i}\in$add$_{\Lambda}U$ for any
$0 \leq i \leq n \}$. We set $U$-resol.dim$_{\Lambda}(A)$ infinity
if no such an integer exists. Dually, for a module $B$ in mod
$\Gamma ^{op}$, we may define $U$-resol.dim$_{\Gamma}(B)$ (see
[2]).

\vspace{0.2cm}

{\bf Lemma 2.7} {\it Let} $E$ {\it be injective in} mod $\Lambda$
({\it resp.} mod $\Gamma ^{op}$). {\it Then} {\it
l.}fd$_{\Gamma}({^*E})$ ({\it resp.} {\it r.}fd$_{\Lambda}({^*E}))
\leq n$ {\it if and only if} $U$-resol.dim$_{\Lambda}(E)$ ({\it
resp.} $U$-resol.dim$_{\Gamma}(E)) \leq n$.

\vspace{0.2cm}

{\it Proof.} Assume that $E$ is injective in mod $\Lambda$ and
{\it l.}fd$_{\Gamma}(^* E)\leq n$. Then there is an exact sequence
$0 \to Q_{n} \to \cdots \to Q_{1} \to Q_{0} \to {^*E} \to 0$ with
each $Q_{i}$ flat (and hence projective) in mod $\Gamma$ for any
$0 \leq i \leq n$. By [3] Chapter VI, Proposition 5.3, ${\rm
Tor}_{j}^{\Gamma}(U , {^*E}) \cong {\rm Hom}_{\Lambda}({\rm
Ext}_{\Gamma}^{j}(U , U), E)=0$ for any $j \geq 1$. Then we easily
have an exact sequence:
$$0 \to U \otimes _{\Gamma}Q_{n} \to \cdots \to
U \otimes _{\Gamma}Q_{1} \to U \otimes _{\Gamma}Q_{0} \to U
\otimes _{\Gamma} {^*E} \to 0.$$ It is clear that $U \otimes
_{\Gamma}Q_{i} \in$add$_{\Lambda}U$ for any $0 \leq i \leq n$. By
[9] p.47, $U \otimes _{\Gamma} {^*E} \cong {\rm
Hom}_{\Lambda}({\rm Hom}_{\Gamma}(U, U),$\linebreak $E) \cong E$.
Hence we conclude that $U$-resol.dim$_{\Lambda}(E)\leq n$.

Conversely, if $U$-resol.dim$_{\Lambda}(E)\leq n$ then there is an
exact sequence $0 \to X_{n} \to \cdots \to X_{1} \to X_{0} \to E
\to 0$ with each $X_{i}$ in add$_{\Lambda}U$ for any $0 \leq i
\leq n$. Since Ext$_{\Lambda}^{j}(U , X_{i})=0$ for any $j \geq 1$
and $0 \leq i \leq n$, $0 \to {^*X}_{n} \to \cdots \to {^*X}_{1}
\to {^*X}_{0} \to {^*E} \to 0$ is exact with each $^*X_{i}$ ($0
\leq i \leq n$) $\Gamma$-projective. Hence we are done.
$\blacksquare$

\vspace{0.2cm}

{\bf Corollary 2.8} {\it Let} $E$ {\it be injective in} mod
$\Lambda$ ({\it resp.} mod $\Gamma ^{op}$). {\it Then} ${^*E}$
{\it is flat in} mod $\Gamma$ ({\it resp.} mod $\Lambda ^{op}$)
{\it if and only if} $_{\Lambda}E\in {\rm add}_{\Lambda}U$ ({\it
resp.} $E_{\Gamma}\in {\rm add}U _{\Gamma})$.

\vspace{0.2cm}

From now on, assume that
$$0 \to {_{\Lambda}U} \buildrel {f_{0}} \over
\longrightarrow E_{0} \buildrel {f_{1}} \over \longrightarrow
E_{1} \buildrel {f_{2}} \over \longrightarrow \cdots \buildrel
{f_{i}} \over \longrightarrow E_{i} \buildrel {f_{i+1}} \over
\longrightarrow \cdots$$ is a minimal injective resolution of
$_{\Lambda}U$.

\vspace{0.2cm}

{\bf Lemma 2.9} {\it Suppose} $U$-dom.dim$(_{\Lambda}U)\geq 1$.
{\it Then, for any} $n \geq 2$, $U$-dom.dim$(_{\Lambda}U )\geq n$
{\it if and only if} grade$_{U}M \geq n$ {\it for any} $M \in$ mod
$\Lambda$ {\it with} $M^*=0$.

\vspace{0.2cm}

{\it Proof.} For any $M\in$mod $\Lambda$ and $i \geq 1$, we have
an exact sequence
$$
{\rm Hom}_{\Lambda}(M, E_{i-1})\to {\rm Hom}_{\Lambda}(M, {\rm
Im}f_{i}) \to {\rm Ext}_{\Lambda}^{i}(M, U) \to 0 \eqno{(\dag)}
$$

Suppose $U$-dom.dim$(_{\Lambda}U )\geq n$. Then $E_{i}$ is
cogenerated by $_{\Lambda}U$ for any $0 \leq i \leq n-1$. So, for
a given $M \in$mod $\Lambda$ with $M^*=0$ we have that
Hom$_{\Lambda}(M, E_{i})=0$ and Hom$_{\Lambda}(M, {\rm
Im}f_{i})=0$ for any $0 \leq i \leq n-1$. Then by the exactness of
(\dag), Ext$_{\Lambda}^{i}(M, U)=0$ for any $1 \leq i \leq n-1$,
and so grade$_{U}M \geq n$.

Now we prove the converse, that is, we will prove that $E_{i}
\in$add$_{\Lambda}U$ for any $0 \leq i \leq n-1$.

First, $E_{0} \in$add$_{\Lambda}U$ by assumption. We next prove
$E_{1} \in$add$_{\Lambda}U$. For any $0\neq x \in$Im$f_{1}$, we
claim that $M^*=$Hom$_{\Lambda}(M, U)\neq 0$, where $M=\Lambda x$.
Otherwise, we have Ext$_{\Lambda}^{i}(M, U)=0$ for any $0 \leq i
\leq n-1$ by assumption. Since $E_{0} \in$add$_{\Lambda}U$,
Hom$_{\Lambda}(M, E_{0})=0$. So from the exactness of (\dag) we
know that Hom$_{\Lambda}(M, {\rm Im}f_{1})=0$, which is a
contradiction. Then we conclude that Im$f_{1}$, and hence $E_{1}$,
is cogenerated by $_{\Lambda}U$. Notice that $E_1$ is finitely
cogenerated,  so $E_{1}\in$add$_{\Lambda}U$. Finally, suppose that
$n \geq 3$ and $E_{i}\in$add$_{\Lambda}U$ for any $0 \leq i \leq
n-2$. Then by using a similar argument to that above we have
$E_{n-1}\in$add$_{\Lambda}U$. The proof is finished.
$\blacksquare$

\vspace{0.2cm}

Dually, we have the following

\vspace{0.2cm}

{\bf Lemma 2.9$'$} {\it Suppose} $U$-dom.dim$(U _{\Gamma})\geq 1$.
{\it Then, for any} $n \geq 2$, $U$-dom.dim$(U _{\Gamma})\geq n$
{\it if and only if} grade$_{U}N \geq n$ {\it for any} $N \in$ mod
$\Gamma ^{op}$ {\it with} $N^*=0$.

\vspace{0.2cm}

We now are in a position to prove the main result in this paper.

\vspace{0.2cm}

{\it Proof of Theorem 1.3.} We only need to prove
$U$-dom.dim$(_{\Lambda}U)\leq U$-dom.dim$(U _{\Gamma})$. Without
loss of generality, suppose $U$-dom.dim$(_{\Lambda}U)=n$.

The case $n=1$ follows from Corollaries 2.5 and 2.8. Let $n \geq
2$. Notice that $U$-dom.dim$(_{\Lambda}U) \geq 1$ and
$U$-dom.dim$(U _{\Gamma}) \geq 1$. By Lemma 2.9$'$ it suffices to
show that grade$_{U}N \geq n$ for any $N \in$ mod $\Gamma ^{op}$
with $N^*=0$. By Lemmas 2.3 and 2.7, for any $i \geq 1$,
Hom$_{\Lambda}$(Ext$_{\Gamma}^{i}(N, U), E_{0}) \cong$
Tor$_{i}^{\Gamma}(N, {^*E}_{0})=0$, so [Ext$_{\Gamma}^{i}(N,
U)]^*=0$. Then by assumption and Lemma 2.9,
grade$_{U}$Ext$_{\Gamma}^{i}(N, U) \geq n$ for any $i \geq 1$. It
follows from Lemma 2.6 that grade$_{U}N \geq n$. $\blacksquare$

\vspace{0.2cm}

\centerline{\bf 3. Some applications}

\vspace{0.2cm}

As applications to the results in above section, we give in this
section some characterizations of $(-)^{**}$ preserving
monomorphisms and being left exact respectively.

\vspace{0.2cm}

Assume that
$$0 \to U _{\Gamma} \buildrel {f_{0}'} \over
\longrightarrow E_{0}' \buildrel {f_{1}'} \over \longrightarrow
E_{1}' \buildrel {f'_{2}} \over \longrightarrow \cdots \buildrel
{f_{i}'} \over \longrightarrow E_{i}' \buildrel {f_{i+1}'} \over
\longrightarrow \cdots$$ is a minimal injective resolution of $U
_{\Gamma}$. We first have the following

\vspace{0.2cm}

{\bf Proposition 3.1} {\it The following statements are equivalent
for any positive integer} $k$.

(1) $U$-dom.dim$(_{\Lambda}U) \geq k$.

(2) $0 \to (_{\Lambda}U)^{**} \buildrel {f_{0}^{**}} \over
\longrightarrow E_{0}^{**} \buildrel {f_{1}^{**}} \over
\longrightarrow E_{1}^{**} \cdots \buildrel {f_{k-1}^{**}} \over
\longrightarrow E_{k-1}^{**}$ {\it is exact}.

(1)$^{op}$ $U$-dom.dim$(U _{\Gamma}) \geq k$.

(2)$^{op}$ $0 \to (U _{\Gamma})^{**} \buildrel {(f'_{0})^{**}}
\over \longrightarrow (E'_{0})^{**} \buildrel {(f'_{1})^{**}}
\over \longrightarrow (E'_{1})^{**} \cdots \buildrel
{(f'_{k-1})^{**}} \over \longrightarrow (E'_{k-1})^{**}$ {\it is
exact}.

\vspace{0.2cm}

{\it Proof.} By Theorem 1.3 we have $(1)\Leftrightarrow (1)^{op}$.
By symmetry, we only need to prove $(1) \Leftrightarrow (2)$.

If $U$-dom.dim$(_{\Lambda}U) \geq k$, then $E_{i}$ is in
add$_{\Lambda}U$ for any $1 \leq i \leq k-1$. Notice that
$_{\Lambda}U$ and each $E_{i}$ ($0 \leq i \leq k-1$) are
$U$-reflexive and hence we have that $0 \to (_{\Lambda}U)^{**}
\buildrel {f_{0}^{**}} \over \longrightarrow E_{0}^{**} \buildrel
{f_{1}^{**}} \over \longrightarrow E_{1}^{**} \cdots \buildrel
{f_{k-1}^{**}} \over \longrightarrow E_{k-1}^{**}$ is exact.
Assume that $(2)$ holds. We proceed by induction on $k$. By
assumption we have the following commutative diagram with exact
rows:
$$\xymatrix{0 \ar[r] & _{\Lambda}U \ar[r]^{f_0} \ar[d]^{\sigma _U} &
E_0 \ar[r]^{f_1} \ar[d]^{\sigma _{E_0}} & E_1 \ar[r]^{f_2}
\ar[d]^{\sigma _{E_1}} &
\cdots \ar[r]^{f_{k-1}} & E_{k-1} \ar[d]^{\sigma _{E_{k-1}}}\\
0 \ar[r] & (_{\Lambda}U)^{**} \ar[r]^{f_0^{**}} & E_0^{**}
\ar[r]^{f_1^{**}} & E_1^{**} \ar[r]^{f_2^{**}} & \cdots
\ar[r]^{f_{k-1}^{**}} & E_{k-1}^{**} }
$$
Since $\sigma _U$ is an isomorphism, $\sigma _{E_{0}}f_0= f_0^{**}
\sigma _U$ is a monomorphism. But $f_0$ is essential, so $\sigma
_{E_{0}}$ is monic, that is, $E_{0}$ is $U$-torsionless and $E_0$
is cogenerated by $_{\Lambda}U$. Moreover, $E_0$ is finitely
cogenerated, so we have that $E_{0} \in$add$_{\Lambda}U$ (and
hence $\sigma _{E_{0}}$ is an isomorphism). The case $k=1$ is
proved. Now suppose that $k\geq 2$ and $E_{i} \in$add$_{\Lambda}U$
(and then $\sigma _{E_{i}}$ is an isomorphism) for any $0\leq i
\leq k-2$ . Put $A_0={_{\Lambda}U}$, $B_0=(_{\Lambda}U)^{**}$,
$g_0=f_0$, $g'_0=f_0^{**}$ and $h_{0}=\sigma _{U}$. Then, for any
$0 \leq i \leq k-2$, we get the following commutative diagrams
with exact rows:
$$\xymatrix{0 \ar[r] & A_i \ar[r]^{g_i} \ar[d]^{h_i} &
E_i \ar[d]^{\sigma _{E_i}} \ar[r] &
A_{i+1} \ar[r] \ar[d]^{h_{i+1}} & 0\\
0 \ar[r] & B_i \ar[r]^{g'_i} &
E_i^{**} \ar[r] & B_{i+1} \ar[r] &
0 }
$$
and
$$\xymatrix{0 \ar[r] & A_{i+1} \ar[r]^{g_{i+1}} \ar[d]^{h_{i+1}} &
E_{i+1} \ar[d]^{\sigma _{E_{i+1}}}\\
0 \ar[r] & B_{i+1} \ar[r]^{g'_{i+1}} & E_{i+1}^{**} }
$$
where $A_{i}={\rm Im}f_{i}$ and $A_{i+1}={\rm Im}f_{i+1}$,
$B_{i}={\rm Im}f_{i}^{**}$ and $B_{i+1}={\rm Im}f_{i+1}^{**}$,
$g_{i}$ and $g_{i+1}$ are essential monomorphisms, $h_{i}$ and
$h_{i+1}$ are induced homomorphisms. We may get inductively that
each $h_j$ is an isomorphism for any $0 \leq j \leq k-1$. Because
$\sigma _{E_{k-1}}g_{k-1}= g'_{k-1} h_{k-1}$ is a monomorphism, by
using a similar argument to that above we have $E_{k-1}
\in$add$_{\Lambda}U$. Hence we conclude that
$U$-dom.dim$(_{\Lambda}U) \geq k$. $\blacksquare$

\vspace{0.2cm}

{\bf Proposition 3.2} {\it The following statements are
equivalent.}

(1) $U$-dom.dim$(_{\Lambda}U) \geq 1$.

(2) $(-)^{**}:$ mod $\Lambda \to$ mod $\Lambda$ {\it preserves
monomorphisms}.

(3) $0 \to (_{\Lambda}U)^{**} \buildrel {f_{0}^{**}} \over
\longrightarrow E_{0}^{**}$ {\it is exact}.

(1)$^{op}$ $U$-dom.dim$(U _{\Gamma}) \geq 1$.

(2)$^{op}$ $(-)^{**}:$ mod $\Gamma ^{op} \to$ mod $\Gamma ^{op}$
{\it preserves monomorphisms}.

(3)$^{op}$ $0 \to (U _{\Gamma})^{**} \buildrel {(f'_{0})^{**}}
\over \longrightarrow (E'_{0})^{**}$ {\it is exact}.

\vspace{0.2cm}

{\it Proof.} By Theorem 1.3 we have $(1)\Leftrightarrow (1)^{op}$.
By symmetry, we only need to prove that the conditions of (1), (2)
and (3) are equivalent.

$(1)\Rightarrow (2)$ If $U$-dom.dim$(_{\Lambda}U) \geq 1$ then
$t(X)=$Ker$\sigma _{X}$ for any $X \in$mod $\Lambda$ by Corollary
2.8 and Proposition 2.4. So $(-)^{**}$ preserves monomorphisms by
Proposition 2.2.

$(2)\Rightarrow (3)$ is trivial and $(3)\Rightarrow (1)$ follows
from Proposition 3.1. $\blacksquare$

\vspace{0.2cm}

{\it Remark.} Proposition 3.2 develops [5] Theorem 1.

\vspace{0.2cm}

{\bf Proposition 3.3} {\it The following statements are
equivalent.}

(1) $U$-dom.dim$(_{\Lambda}U) \geq 2$.

(2) $(-)^{**}:$ mod $\Lambda \to$ mod $\Lambda$ {\it is left
exact}.

(3) $0 \to (_{\Lambda}U)^{**} \buildrel {f_{0}^{**}} \over
\longrightarrow E_{0}^{**} \buildrel {f_{1}^{**}} \over
\longrightarrow E_{1}^{**}$ {\it is exact}.

(4) $(-)^{**}:$ mod $\Lambda \to$ mod $\Lambda$ {\it preserves
monomorphisms and} Ext$_{\Gamma}^{1}({\rm Ext} _{\Lambda}^{1}(X,
U), U)=0$ {\it for any} $X\in$mod $\Lambda$.

(1)$^{op}$ $U$-dom.dim$(U _{\Gamma}) \geq 2$.

(2)$^{op}$ $(-)^{**}:$ mod $\Gamma ^{op} \to$ mod $\Gamma ^{op}$
{\it is left exact}.

(3)$^{op}$ $0 \to (U _{\Gamma})^{**} \buildrel {(f'_{0})^{**}}
\over \longrightarrow (E'_{0})^{**} \buildrel {(f'_{1})^{**}}
\over \longrightarrow (E'_{1})^{**}$ {\it is exact}.

(4)$^{op}$ $(-)^{**}:$ mod $\Gamma ^{op}\to$ mod $\Gamma ^{op}$
{\it preserves monomorphisms and} Ext$_{\Lambda}^{1}({\rm Ext}
_{\Gamma}^{1}(Y, U), U)=0$ {\it for any} $Y\in$mod $\Gamma ^{op}$.

\vspace{0.2cm}

{\it Proof.} By Theorem 1.3 we have $(1)\Leftrightarrow (1)^{op}$
and by Proposition 3.1 we have $(1)\Leftrightarrow (3)$. So, by
symmetry we only need to prove that $(1)\Leftrightarrow (2)$ and
$(1)\Rightarrow (4) \Rightarrow (1)^{op}$.

$(1) \Leftrightarrow (2)$ Assume that $(-)^{**}:$ mod $\Lambda
\to$ mod $\Lambda$ is left exact. Then, by Proposition 3.2, we
have that $U$-dom.dim$(_{\Lambda}U) \geq 1$ and $E_{0}
\in$add$_{\Lambda}U$.

Let $K=$Im$(E_{0} \to E_{1})$ and $v: K \to E_{1}$ be the
essential monomorphism. By assumption and the exactness of the
sequences $0 \to U \to E_{0} \to K \to 0$ and $0 \to K \buildrel
{v} \over \longrightarrow E_{1}$, we have the following exact
commutative diagrams:

$$\xymatrix{0 \ar[r] & U \ar[r] \ar[d]^{\sigma _U} &
E_0 \ar[r] \ar[d]^{\sigma _{E_0}} &
K \ar[r] \ar[d]^{\sigma _K} & 0\\
0 \ar[r] & U^{**} \ar[r] & E_0^{**} \ar[r] & K^{**} }
$$
and

$$\xymatrix{0 \ar[r] & K \ar[r]^v \ar[d]^{\sigma _K} &
E_1 \ar[d]^{\sigma _{E_1}}\\
0 \ar[r] & K^{**} \ar[r]^{v^{**}} & E_1^{**} }
$$
where $\sigma _{U}$ and $\sigma _{E_{0}}$ are isomorphisms. By
applying the snake lemma to the first diagram we have that $\sigma
_{K}$ is monic. Then we know from the second diagram that $\sigma
_{E_{1}}v= v^{**} \sigma _{K}$ is a monomorphism. However, $v$ is
essential, so $\sigma _{E_{1}}$ is monic, that is, $E_{1}$ is
$U$-torsionless and $E_1$ is cogenerated by $_{\Lambda}U$.
Moreover, $E_1$ is finitely cogenerated, so we conclude that
$E_{1} \in$add$_{\Lambda}U$.

Conversely, assume that $U$-dom.dim$(_{\Lambda}U) \geq 2$ and $0
\to A \buildrel {\alpha} \over \longrightarrow B \buildrel {\beta}
\over \longrightarrow C \to 0$ is an exact sequence in mod
$\Lambda$. By Proposition 3.2, $\alpha ^{**}$ is monic. By
assumption, Corollary 2.8 and Lemma 2.3 we have Hom$_{\Gamma}({\rm
Ext}_{\Lambda}^{1}(C, U), E_{0})=0$. Since Coker$\alpha ^*$ is
isomorphic to a submodule of Ext$_{\Lambda}^{1}(C, U)$,
Hom$_{\Gamma} ({\rm Coker}\alpha ^*, E_{0})=0$ and
Hom$_{\Gamma}({\rm Coker}\alpha ^*, U)=0$. Then, by Theorem 1.3
and Lemma 2.9$'$, grade$_{U}$Coker$\alpha ^* \geq 2$. It follows
easily that $0 \to A^{**} \buildrel {\alpha ^{**}} \over
\longrightarrow B^{**} \buildrel {\beta ^{**}} \over
\longrightarrow C^{**}$ is exact.

$(1) \Rightarrow (4)$ Suppose $U$-dom.dim$(_{\Lambda}U) \geq 2$.
By Proposition 3.2, $(-)^{**}:$ mod $\Lambda \to$ mod $\Lambda$
preserves monomorphisms. On the other hand, we have that
$U$-dom.dim$(U _{\Gamma}) \geq 2$ by Theorem 1.3. It follows from
Corollary 2.8 and Lemma 2.3 that Hom$_{\Gamma}({\rm
Ext}_{\Lambda}^{1}(X, U), E_{0}^{'})=0$ for any $X\in$mod
$\Lambda$. So $[{\rm Ext}_{\Lambda}^{1}(X, U)]^*=0$ and hence
Ext$_{\Gamma}^{1}({\rm Ext} _{\Lambda}^{1}(X, U), U)=0$ by Lemma
2.9$'$.

$(4)\Rightarrow (1)^{op}$ Suppose that (4) holds. Then
$U$-dom.dim$(U _{\Gamma}) \geq 1$ by Proposition 3.2.

Let $A$ be in mod $\Lambda$ and $B$ any submodule of
Ext$_{\Lambda}^{1}(A, U)$ in mod $\Gamma ^{op}$. Since
$U$-dom.dim$(U _{\Gamma})$\linebreak $\geq 1$, Hom$_{\Gamma}({\rm
Ext}_{\Lambda}^{1}(A, U), E_{0}^{'})=0$ by Corollary 2.8 and Lemma
2.3. So Hom$_{\Gamma}(B, E_{0}^{'})=0$ and hence Hom$_{\Gamma}(B,
E_{0}^{'}/U _{\Gamma})\cong {\rm Ext}_{\Gamma}^{1}(B, U
_{\Gamma})$. On the other hand, Hom$_{\Gamma}(B, E_{0}^{'})=0$
implies $B^*=0$. Then by [6] Lemma 2.1 we have that
$B\cong$Ext$_{\Lambda}^{1}({\rm Tr}_{U}B, U)$ with ${\rm Tr}_{U}B$
in mod $\Lambda$. By (4), Hom$_{\Gamma}(B, E_{0}^{'}/U)\cong {\rm
Ext}_{\Gamma}^{1}(B, U) \cong {\rm Ext}_{\Gamma}^{1}({\rm
Ext}_{\Lambda}^{1}({\rm Tr}_{U}B, U), U)=0$. Then by using a
similar argument to that in the proof $(2)^{op}\Rightarrow
(1)^{op} $ in Proposition 2.2, we have that Hom$_{\Gamma}({\rm
Ext}_{\Lambda}^{1}(A, U), E_{1}^{'})=0$ (note: $E_{1}^{'}$ is the
injective envelope of $E_{0}^{'}/U$). Thus $E_{1}^{'}\in$add$U
_{\Gamma}$ by Lemma 2.3 and Corollary 2.8, and therefore
$U$-dom.dim$(U _{\Gamma}) \geq 2$. $\blacksquare$

\vspace{0.2cm}

{\it Remark.} Proposition 3.3 develops [5] Theorem 2.

\vspace{0.2cm}

Finally we give some equivalent characterizations of
$U$-resol.dim$_{\Lambda}(E_{0}) \leq 1$ as follows.

\vspace{0.2cm}

{\bf Proposition 3.4} {\it The following statements are
equivalent.}

(1) $U$-resol.dim$_{\Lambda}(E_{0}) \leq 1$.

(2) $\sigma _{X}$ {\it is an essential monomorphism for any}
$U$-{\it torsionless module} $X$ {\it in} mod $\Lambda$.

(3) $f^{**}$ {\it is a monomorphism for any monomorphism} $f:X\to
Y$ {\it in} mod $\Lambda$ {\it with} $Y$ $U$-{\it torsionless}.

(4) grade$_{U}$Ext$^{1}_{\Lambda}(X, U) \geq 1$ ({\it that is},
$[{\rm Ext}^{1}_{\Lambda}(X, U)]^*=0$) {\it for any} $X$ {\it in}
mod $\Lambda$.

\vspace{0.2cm}

{\it Proof.} $(1) \Rightarrow (2)$ Assume that $X$ is
$U$-torsionless in mod $\Lambda$. Then Coker$\sigma
_{X}\cong$\linebreak ${\rm Ext}_{\Gamma}^{2}({\rm Tr}_{U}X, U)$ by
[6] Lemma 2.1. By Lemmas 2.7 and 2.3 we have
Hom$_{\Lambda}($Coker$\sigma _{X}, E_{0})=$
Hom$_{\Lambda}($Ext$_{\Gamma}^{2}($Tr$_{U}X, U), E_{0})=0$. Then
Hom$_{\Lambda}(A, {_{\Lambda}U})=0$ for any submodule $A$ of
Coker$\sigma _{X}$, which implies that any non-zero submodule of
Coker$\sigma _{X}$ is not $U$-torsionless.

Let $B$ be a submodule of $X^{**}$ with $X \bigcap B=0$. Then $B
\cong B/X\bigcap B \cong (X+B)/X$ is isomorphic to a submodule of
Coker$\sigma _{X}$. On the other hand, $B$ is clearly
$U$-torsionless. So $B=0$ and hence $\sigma _{X}$ is essential.

$(2) \Rightarrow (3)$ Let $f: X\to Y$ be monic in mod $\Lambda$
with $Y$ $U$-torsionless. Then $f^{**}\sigma _{X}=\sigma _{Y}f$ is
monic. By (2), $\sigma _{X}$ is an essential monomorphism, so
$f^{**}$ is monic.

$(3) \Rightarrow (4)$ Let $X$ be in mod $\Lambda$ and $0 \to Y
\buildrel {g} \over \longrightarrow P \to X \to 0$ an exact
sequence in mod $\Lambda$ with $P$ projective. It is easy to see
that $[{\rm Ext}^{1}_{\Lambda}(X, U)]^* \cong {\rm Ker}g^{**}$. On
the other hand, $g^{**}$ is monic by (3). So Ker$g^{**}=0$ and
$[{\rm Ext}^{1}_{\Lambda}(X, U)]^*=0$.

$(4) \Rightarrow (1)$ Let $M$ be in mod $\Gamma ^{op}$ and $\cdots
\to P_{1} \to P_{0} \to M \to 0$ a projective resolution of $M$ in
mod $\Gamma ^{op}$. Put $N=$Coker$(P_{2}\to P_{1})$. By [6] Lemma
2.1, ${\rm Ext}^{2}_{\Gamma}(M, U) \cong {\rm Ext}^{1}_{\Gamma}(N,
U) \cong {\rm Ker}\sigma _{{\rm Tr}_{U}N}$. On the other hand,
since $N$ is $U$-torsionless, ${\rm Ext}^{1}_{\Lambda}({\rm
Tr}_{U}N, U)\cong {\rm Ker}\sigma _{N}=0$.

Let $X$ be any finitely generated submodule of ${\rm
Ext}^{2}_{\Gamma}(M, U)$ and $f_{1}: X\to{\rm Ext}^{2}_{\Gamma}(M,
U) (\cong {\rm Ker}\sigma _{{\rm Tr}_{U}N})$ the inclusion, and
let $f$ be the composition: $X \buildrel {f_{1}} \over
\longrightarrow {\rm Ext}^{2}_{\Gamma}(M, U) \buildrel {g} \over
\longrightarrow {\rm Tr}_{U}N$, where $g$ is a monomorphism. By
using the same argument as that in the proof of $(2)^{op}
\Rightarrow (1)^{op}$ in Proposition 2.2, we get that $f^*=0$.
Hence, by applying Hom$_{\Lambda}(-, U)$ to the exact sequence $0
\to X \buildrel {f} \over \longrightarrow {\rm Tr}_{U}N \to {\rm
Coker}f \to 0$, we have $X^* \cong {\rm Ext}^{1}_{\Lambda}({\rm
Coker}f, U)$. Then $X^{**} \cong [{\rm Ext}^{1}_{\Lambda}({\rm
Coker}f, U)]^*=0$ by (4), which implies that $X^*=0$ since $X^*$
is a direct summand of $X^{***}(=0)$ by [1] Proposition 20.24.
Also by using the same argument as that in the proof of $(2)^{op}
\Rightarrow (1)^{op}$ in Proposition 2.2, we get that
Hom$_{\Lambda}($Ext$^{2}_{\Gamma}(M, U), E_{0})=0$. It follows
from Lemma 2.3 that {\it l.}fd$_{\Gamma}({^*E}_{0})\leq 1$.
Therefore $U$-resol.dim$_{\Lambda}(E_{0}) \leq 1$ by Lemma 2.7.
$\blacksquare$

\vspace{0.2cm}

{\it Remark.} By Theorem 1.3, we have that $E_0
\in$add$_{\Lambda}U$ if and only if $E'_0 \in$add$U_{\Gamma}$,
that is, $U$-resol.dim$_{\Lambda}(E_{0})=0$ if and only if
$U$-resol.dim$_{\Gamma}(E'_{0})=0$. However, in general, we don't
have the fact that $U$-resol.dim$_{\Lambda}(E_{0}) \leq 1$ if and
only if $U$-resol.dim$_{\Gamma}(E'_{0}) \leq 1$ even when
$_{\Lambda}U_{\Gamma}={_{\Lambda}\Lambda _{\Lambda}}$. We use
$I_0$ and $I'_0$ to denote the injective envelope of
$_{\Lambda}\Lambda$ and $\Lambda _{\Lambda}$, respectively.
Consider the following example. Let $K$ be a field and $\Delta$
the quiver:
$$\xymatrix{ 1 \ar @<2pt> [r]^{\alpha} &
2 \ar @<2pt> [l]^{\beta} \ar[r]^{\gamma} & 3 }
$$
(1) If $\Lambda =K\Delta /(\alpha \beta \alpha)$. Then {\it
l.}fd$_{\Lambda}(I_0) =1$ and {\it r.}fd$_{\Lambda}(I'_0) \geq 2$.
(2) If $\Lambda =K\Delta /(\gamma \alpha , \beta \alpha)$. Then
{\it l.}fd$_{\Lambda}(I_0) =2$ and {\it r.}fd$_{\Lambda}(I'_0)
=1$.

\vspace{0.5cm}

{\bf Acknowledgements} The author thanks Prof. Kent. R. Fuller and
the referee for their helpful comments. This work was partially
supported by National Natural Science Foundation of China (Grant
No. 10001017).

\end{document}